\documentclass[a4paper,12pt]{article}
\usepackage{t1enc}
\usepackage[latin2]{inputenc}
\usepackage[T1]{fontenc}
\usepackage{amsmath}
\usepackage{amssymb}
\usepackage{amsthm}
\usepackage{graphicx}
\newtheorem{thm}{Theorem}
\newtheorem{prop}{Proposition}

\newtheorem*{remark}{Remark}

\newcommand{\Z}{\mathbb Z}
\newcommand{\R}{\mathbb R}
\renewcommand{\P}{\mathbf P}
\newcommand{\E}{\mathbf E}
\newcommand{\I}{1\!\!1}

\author{Istv\'an R\'edl \and B\'alint Vet\H o}

\title{Random Walk in Periodic Environment}

\begin{document}

\maketitle

\begin{abstract}

We consider a special case of random walk in random environment (RWRE) on
$\Z^d$ where the environment is periodic (RWPE). Under natural conditions, we
show that law of large numbers and central limit theorem holds. In the
ballistic nearest neighbour reversible case, we prove that the angle between
the asymptotic direction of the RWPE and the average negative gradient of the
potential function of the reversible environment is less than $\pi/2$, that is,
the potential cannot increase asymptotically along the trajectory of the RWPE.
But this angle can be close to $\pi/2$.

\end{abstract}

\section{Introduction}

The general random walk in random environment (RWRE) on $\Z^d$ consists of two
components: a random environment is chosen and a particle performs a random
walk with transition probabilities given by the environment. More precisely,
let $\mathcal M$ be the set of probability measures on $\Z^d$. The environment
is a random element of $\mathcal M^{\Z^d}$, i.e.\ for each $x\in\Z^d$, a random
probability measure $(p_x(y))_{y\in\Z^d}$. Conditionally given the environment,
the random walk $X_n$ is a Markov chain, that is
\[\P\left(X_{n+1}=x+y\bigm|X_n=x,(p_x(z))_{z\in\Z^d}\right)=p_x(y)\]
and $X_0=0$. For more about RWRE, see the lecture notes \cite{OfZe04}.
Motivated partly by \cite{NB}, we describe here the possible asymptotic
directions in a special model for RWRE.

In the present paper, our main assumption is that the environment is periodic:
there are fixed integers $M_1,M_2,\dots,M_d$ such that
\begin{equation}
p_{x+M_ie_i}(y)=p_x(y)\qquad i=1,\dots,d\label{periodicity}
\end{equation}
for all $x,y\in\Z^d$ where $e_i$ is the $i$th unit vector in $\Z^d$. We
consider a fixed frozen environment, because our theorems are true for almost
all realizations of the environment.

With the notation
\[M:=M_1\Z\times\dots\times M_d\Z,\]
we can introduce the following equivalence relation of vertices in $\Z^d$:
\[x\sim y\quad\mbox{if}\quad x-y\in M.\]
Condition \eqref{periodicity} is equivalent with saying that, for equivalent
vertices, the step distributions are the same.

Let
\[T:=\Z^d/M=\{0,\dots,M_1-1\}\times\dots\times\{0,\dots,M_d-1\}\]
be the torus obtained by factorizing the lattice with the sublattice $M$. In
this way, we get a finite set. Note that, by the periodicity condition, the
measures
\[(p_x(y))_{y\in\Z^d}\qquad x\in T\]
extend to an environment on the whole $\Z^d$.

Furthermore, let $\xi_n$ be the \emph{induced Markov chain} on the finite state
space $T$ defined by taking the equivalence class of $X_n$. We call the
transition matrix of $\xi_n$ by $P$, i.e.
\[P_{xy}:=\P\left(\xi_n=y\bigm|\xi_{n-1}=x\right)\]
for any $x,y\in T$. We suppose that the environment is given in such a way that
$\xi_n$ is irreducible and aperiodic, hence, it has a unique stationary
probability measure $\pi$ on $T$.

In this setting, the increments $(X_n-X_{n-1})_{n\ge1}$ form a sequence of
\emph{chain dependent random variables}. They actually depend on the pair
$(\xi_{n-1},\xi_n)$, which is also a finite state Markov chain. Indeed, given
that the $n$th jump $X_{n-1}\mapsto X_n$ starts at the equivalence class
$\xi_{n-1}$ and that it ends in the equivalence class $\xi_n$, then the
distribution of $X_n-X_{n-1}$ is
\begin{equation}
\P\left(X_n-X_{n-1}=w\bigm|\xi_{n-1}=x,\xi_n=y\right) =\I(w\sim
y)\frac{p_x(w)}{\sum_{z\in\Z^d:z\sim y} p_x(z)}\label{conddistr}
\end{equation}
by taking conditional probability only for those pairs $(x,y)$ for which the
$x\mapsto y$ jump has positive probability. Therefore, the sequence $X_n$ can
be sampled as follows. One generates first the whole trajectory of the induced
Markov chain $(\xi_n)$, then, for each $n$, conditionally on the pair
$(\xi_{n-1},\xi_n)$, one draws $X_n-X_{n-1}$ independently with the appropriate
distribution given by \eqref{conddistr}. For any two states $x,y\in T$, let
$\mu_{xy}$ be the conditional expectation of the vector of the $x\mapsto y$
jump for the RWPE $X_n$, that is, the expectation of \eqref{conddistr}:
\[\mu_{xy}:=\E\left(X_n-X_{n-1}\bigm|\xi_{n-1}=x,\xi_n=y\right)
=\frac{\sum_{z\in\Z^d:z\sim y} zp_x(z)}{\sum_{z\in\Z^d:z\sim y} p_x(z)}.\]
If no
$x\mapsto y$ jump is possible, then the conditional distribution does not make
sense, and let $\mu_{xy}=0$.

For chain dependent random variables, see the textbook \cite{Durr96}. They are
also referred to as random walks with internal degrees of freedom, see
\cite{SzaszD83}.

This paper is organized as follows: in Section \ref{sec:res}, we state the
results: under natural moment conditions, law of large numbers and central
limit theorem hold in this setup. For the reversible nearest neighbour RWPE in
the ballistic case, we investigate the relationship between the average
gradient of the potential and the asymptotic direction of the random walk. The
proofs are done in Section \ref{sec:lln,clt} and \ref{sec:rev}.

\section{Results}
\label{sec:res}

\begin{prop}[Law of large numbers]\label{prop:lln}
Let $X_n$ be a RWPE. Suppose that the induced Markov chain $\xi_n$ is
irreducible and aperiodic, and, for the expectations,
\[\sum_{y\in\Z^d}|y|p_x(y)<\infty\]
holds for all $x\in T$. Then
\[\frac{X_n}n\to\nu\]
almost surely as $n\to\infty$ where
\begin{equation}
\nu:=\sum_{x\in T} \pi(x)\sum_{y\in\Z^d} p_x(y)y=\sum_{x,y\in T}
\pi(x)P_{xy}\mu_{xy}.\label{defnu}
\end{equation}
\end{prop}

\begin{prop}[Central limit theorem]\label{prop:clt}
Let $X_n$ be a RWPE. Suppose that the induced Markov chain $\xi_n$ is
irreducible and aperiodic, and, for the second moments,
\[\sum_{y\in\Z^d}|y|^2p_x(y)<\infty\]
for all $x\in T$. Then
\[\frac{X_n-n\nu}{\sqrt n}\Longrightarrow\mathcal N(0,\Sigma)\]
in distribution as $n\to\infty$ where
\begin{align}
\Sigma_{ij}&:=\sum_{x\in T} \pi(x) \sum_{y\in\Z^d} p_x(y)
(y_i-\nu_i)(y_j-\nu_j)\label{defSigma}\\
&\quad+2\sum_{x,y,z,w\in T} \pi(x)P_{xy}((\mu_{xy})_i-\nu_i) (I-P)^{-1}_{yz} P_{zw}
((\mu_{zw})_j-\nu_j).\notag
\end{align}
\end{prop}

These theorems are consequences of general theory of Markov chians. For the
central limit theorem, our proof is remarkably simpler than the one found in
\cite{Take02}.

The main result of this paper is about nearest neighbour RWPE which is
reversible with respect to some $\sigma$-finite measure. The reversibility
condition is satisfied if, for each elementary cirle of the lattice, the
products of the jump probabilities in both directions are equal.

The nearest neighbour property means that only the probabilities $p_x(e)$ are
non-zero where $|e|=1$. We assume also that these probabilities are strictly
positive, hence, irreducibility of $\xi_n$ follows. However, we lose the
aperiodicity assumption on $\xi_n$.

Reversibility means that we can define the potential function of the
environment by $u(0)=0$ and
\[u(x)-u(x+e)=\log\left(\frac{p_x(e)}{p_{x+e}(-e)}\right)\]
for each $x\in\Z^d$ and $|e|=1$. The \emph{average negative gradient} of the
potential function $u$ is given by
\[g:=\left(\frac1{M_1}\log \prod_{i=0}^{M_1-1} \frac{p_{ie_1}(e_1)}{p_{(i+1)e_1}(-e_1)},\dots,
\frac1{M_d}\log\prod_{i=0}^{M_d-1}\frac{p_{ie_d}(e_d)}{p_{(i+1)e_d}(-e_d)}\right).\]

\begin{thm}\label{thm:main}
Let $X_n$ be a reversible nearest neighbour RWPE. Assume that $p_x(e)>0$ for
all $x\in\Z^d$ and $|e|=1$. If $g$ is not the $0$ vector, then $\langle
g,\nu\rangle>0$ where $\nu$ is given by \eqref{defnu}.
\end{thm}

This theorem tells us that the potential of the random walker can only decrease
asymptotically.

\begin{remark}
Theorem \ref{thm:main} is the best we can prove, because the angle of $g$ and
$\nu$ can be arbitrarily close to $\pi/2$ even in $\Z^2$ in the simple random
walk case. Consider the jump distribution

\begin{align*}
p_x(e_1)&=K\varepsilon, & p_x(-e_1)&=\varepsilon,\\
p_x(e_2)&=\frac23(1-(K+1)\varepsilon), & p_x(-e_2)&=\frac13(1-(K+1)\varepsilon)
\end{align*}
for all $x\in\Z^2$ with $K$ and $\varepsilon$ specified later.

Then the (average) negative gradient vector of the potential is clearly $(\log
K,\log2)$ for any $\varepsilon$ as long as all the weights are positive. So let
us choose $K$ first. We can make the angle of the gradient $(\log K,\log2)$
with the first coordinate vector $e_1$ arbitrarily small.

On the other hand, the expected jump vector is obviously
\[\left((K-1)\varepsilon,\frac13(1-(K+1)\varepsilon)\right)
\to\left(0,\frac13\right)\quad\mbox{as}\quad\varepsilon\to0\]
for any fixed $K$.
\end{remark}

The proof of Theorem \ref{thm:main} is based on the following idea: we
introduce level hyperplanes which are perpendicular to the gradient vector $g$.
The RWPE taken in the hitting times of these level hyperplanes is a one
dimensional simple random walk, which is transient in the asymmetric case.

\begin{remark}
The one dimensional case of the RWPE can be also handled by standard martingale
technics. The key observation here is that the RWPE stopped at the exit time of
some interval $[-K,K]$ is a generalized gamblers ruin game.
\end{remark}

\section{Law of large numbers and central limit theorem}
\label{sec:lln,clt}

The proof of Proposition \ref{prop:lln} is easy: $\xi_n$ converges
exponentially fast to its stationary distribution. Conditionally given the
trajectory of $\xi_n$, for any $x,y\in T$, those jumps of $X_n$ which start in
the equivalence class $x$ and end in the equivalence class $y$ are i.i.d.\
random variables with mean $\mu_{xy}$. The asymptotic proportion of these jumps
is close to the stationary weight of the $x\mapsto y$ jumps of $\xi_n$. Further
details are left to the reader.

Proposition \ref{prop:clt} follows from Exercise 7.3 in Chapter 7 of
\cite{Durr96}. It gives the central limit theorem for general sequences of
chain dependent random variables. Our setup fits into this framework. The
covariance matrix is given by the discrete version of the Green\,--\,Kubo
formula
\[\E(X_1-X_0-\nu)^2+2\sum_{n=1}^\infty \E(X_{n+1}-X_n-\nu)(X_1-X_0-\nu),\]
which transforms to \eqref{defSigma} in our case.

\section{Asymptotic direction and average gradient of the environment}
\label{sec:rev}

We say that the vector $g$ is \emph{appropriate}, if it has rational
coordinates and it has a multiple $g_1=\alpha g$ with $\alpha\in\R$ such that
the endpoint of $g_1$ is equivalent to its starting point. Note that
appropriate directions are dense in $\R^d$. We assume first that $g$
appropriate. Later this assumption will be relaxed.

We take the hyperplane across the origin orthogonal to the gradient $g$, that
is, $H_0:=\{x\in\R^d:\langle x,g \rangle=0\}$. Let $g_1$ be the vector given
above. We consider
\[H_k:=\{x\in\R^d:\langle x,g \rangle=k\langle g_1,g\rangle\}=H_0+kg_1\qquad
k\in\Z\]
the hyperplane $H_0$ translated by the multiples of $g_1$.

We introduce a new equivalence relation which is a refinement of $\sim$, i.e.\
let
\[x\approx y\qquad\mbox{if and only if}\qquad x\sim y\quad\mbox{and}\quad
u(x)=u(y),\]
that is, if $x$ and $y$ are $\sim$ equivalent and they also have the same
potential value. Note that the classes are non-trivial, since $g$ has rational
coordinates.

We can now define the \emph{level hyperplane}s $L_k$ in $\Z^d$ which are the
discrete counterparts of $H_k$. Let $L_0$ be a subset in $\Z^d$ with the
following properties:
\begin{enumerate}
\item $L_0$ is periodic, i.e.\ it contains entire equivalence classes of
$\approx$;
\item the distances of its points from $H_0$ is uniformly bounded;
\item $L_0$ disconnects $\Z^d$ into two components.
\end{enumerate}
It is not hard to construct such a set. The level hyperplanes are given by
\begin{equation}
L_k:=L_0+kg_1.\label{defLk}
\end{equation}
Note that, by the second property, finitely many $\approx$ equivalence classes
are in $L_k$ for each $k$, therefore, the number of possible potential values
in $L_k$ is finite. This is again a consequence of the rational coordinates in
$g$.

For each pair of $\approx$ equivalence classes that both appear in $L_0$, we
choose a directed path $\gamma_{x,y}$ where $x$ is the starting point, $y$ is
the endpoint, which are in the given two classes. This construction gives us
finitely many paths, hence their weight is uniformly bounded from below:
\begin{equation}
\P(\gamma_{x,y})\ge\beta>0\qquad x,y\in L_0.\label{gammabound}
\end{equation}
Moreover, by the definition \eqref{defLk}, each pair of $\approx$ equivalence
classes in any of the $L_k$s can be connected by a shifted version of some
$\gamma_{x,y}$ with $x,y\in L_0$. In this way, the $\gamma_{x,y}$ is extended
to $x,y\in L_k$ for all $k$ and the boundedness property \eqref{gammabound} is
preserved.

Let
\[\tau_k:=\min\{n\ge0:X_n\in L_k\}\]
be the first hitting time of $L_k$. We will prove that, for any $x_0\in L_0$,
\begin{equation}
\P\left(\tau_{-k}<\tau_k\bigm|X_0=x_0\right)<\frac12\label{hittingprob}
\end{equation}
for some $k>0$ integer. If this relation is verified, then the RWPE is an
asymmetric simple random walk on the level hyperplanes $(L_{kl})_{l\in\Z}$
where the time needed for one step ($\tau_k\wedge\tau_{-k}$) can be dominated
by a geometrically distributed random variable, and the periodicity provides
independence between two neighbouring levels. This is enough for the proof of
Theorem \ref{thm:main}.

For showing \eqref{hittingprob}, we give a uniform estimate on the probability
of paths going from $L_0$ to $L_{-k}$ where $k>0$ is some integer for which
$L_0$ and $L_{-k}$ are disjoint. We take an arbitrary path from $L_0$ to
$L_{-k}$. We call its starting point $a$, its last point in $L_0$ is $b$ ($b=a$
is also possible), and the endpoint of the path is $c$, which is also its first
point in $L_{-k}$. We denote the path by
\[\omega_{a,c}=\omega_{a,b}*\omega_{b,c}\]
where $*$ means connection of two paths.

We construct a path from $L_0$ to $L_k$ that serves for estimation. The new
path starts from $a\in L_0$. The first part is one of the previously fixed
paths $\gamma_{a,b'}$ where $b'\in L_0$ with $b'\sim b$. Then, we take the
shifted and reversed version of the path $\omega_{a,b}$ which we call
$\upsilon_{b',a'}$ where $a'\in L_0$ and $a'\sim a$. In the next step, we
continue along the canonically fixed path $\gamma_{a',c'}$ where $c'\in L_0$
and $c'\sim c$. Finally, we shift and reverse the path $\omega_{b,c}$ to get
$\upsilon_{c',b''}$ where $b''\in L_k$ with $b''\sim b$. With this procedure,
we have constructed the path
\[\psi_{a,b''}:=\gamma_{a,b'}*\upsilon_{b',a'}*\gamma_{a',c'}*\upsilon_{c',b''}\]
which starts at $L_0$ and it arrives at $L_k$. Furthermore, no other points of
$\psi_{a,b''}$ than $b''$ is contained in $L_k$, recall the definition of the
point $b$.

By the definition of the potential and the periodicity, it follows that
\[\P(\psi_{a,b''})
=e^{u(c)-u(a)}\P(\omega_{a,c})\P(\gamma_{a,b'})\P(\gamma_{a',c'}).\]
From
\eqref{gammabound},
\begin{equation}
\P(\omega_{a,c}) \le
e^{u(a)-u(c)}\beta^{-2}\P(\psi_{a,b''}).\label{weigthestimate}
\end{equation}

It is not hard to show that the map $\omega_{a,c}\mapsto\psi_{a,b''}$ is
injective. If we choose $k$ large enough, then the factor
$e^{u(a)-u(c)}\beta^{-2}$ in \eqref{weigthestimate} can be made smaller than
$1/2$ for all pairs $a\in L_0$ and $c\in L_k$. Therefore, all paths
$\omega_{a,c}$ is associated injectively to a path $\psi_{a,b''}$ with strictly
higher probability by \eqref{weigthestimate}. This proves \eqref{hittingprob}.

If $g$ is not appropriate, then we approximate $g$ with some appropriate $g'$.
We build up the level hyperplanes using $g'$. The error can be made
arbitrarily small, because the walker cannot get too far in other directions
before hitting $L_{-k}$ by large deviation principle. This completes the
proof of Theorem \ref{thm:main}.\\

\noindent {\bf Acknowledgement:} We thank Bal\'azs R\'ath for drawing our
attention to this problem.


\begin{thebibliography}{9}

\bibitem{NB}
N.\ Berger: Slowdown estimates for ballistic random walk in random environment.
{\it Preprint}, \texttt{http://arxiv.org/abs/0811.1710} (2010)

\bibitem{Durr96}
R.\ Durrett: {\it Probability Theory and Examples}. Second Ed., Duxbury Press,
1996.

\bibitem{SzaszD83}
A.\ Krámli, D.\ Szász: Random Walks with Internal Degrees of Freedom. {\it
Zeitschrift für Wahrscheinlichkeitstheorie und verwandte Gebiete}, {\bf
63}:85--95 (1983)

\bibitem{Take02}
T.\ Takenami: Local limit theorem for random walk in periodic environment. {\it
Osaka J.\ Math.} \textbf{39}:867--895 (2002)

\bibitem{OfZe04}
O.\ Zeitouni: {\it Lecture Notes on Random Walks in Random Environment.}
Lecture Notes in Mathematics, vol.\ 1837 pp.\ 190--312 Springer, Heidelberg,
2004.

\end{thebibliography}
\end{document}